\documentclass[12pt]{article}
\usepackage{mathrsfs}
\usepackage{stmaryrd}
\usepackage{amsfonts,amsmath,amssymb,amscd}
\usepackage{shadow}
\usepackage{graphicx,floatrow}
\usepackage{color}
\usepackage{longtable}
\usepackage{pstricks,multido}
\usepackage{hyperref}
\usepackage{qtree}
\usepackage{tabularx}
\usepackage{array}
\usepackage{caption}
\allowdisplaybreaks

\newtheorem{thm}{Theorem}[section]
\newtheorem{lem}[thm]{Lemma}

\newtheorem{cor}[thm]{Corollary}
\newtheorem{conj}[thm]{Conjecture}

\def\qed{\hfill \rule{4pt}{7pt}}

\def\pf{\noindent {\it{Proof.} \hskip 2pt}}

\numberwithin{equation}{section}

\begin{document}
\begin{center}
{\large\bf   Equidistributions of MAJ and STAT over
pattern avoiding permutations}
\end{center}

\begin{center}
Joanna N. Chen

College of Science\\
Tianjin University of Technology\\
Tianjin 300384, P.R. China

joannachen@tjut.edu.cn.

\end{center}

\begin{abstract}
Babson and Steingr\'{\i}msson introduced generalized permutation patterns and showed that most of the Mahonian statistics in the literature can be expressed by
the combination of generalized pattern functions. Particularly, they defined
a new Mahonian statistic in terms of generalized pattern functions, which is denoted  $stat$. Recently,
Amini investigated the equidistributions of  these Mahonian
statistics over sets of pattern avoiding permutations.
Moreover, he posed  several conjectures. In this paper, we construct
a bijection from $S_n(213)$ to $S_n(231)$, which maps
the statistic $(maj,stat)$ to the statistic $(stat,maj)$. This allows us to give  solutions to some of Amini's conjectures.
\end{abstract}

\noindent {\bf Keywords}: Mahonian statistic, pattern, generalize pattern, bijection

\noindent {\bf AMS  Subject Classifications}: 05A05, 05A15

\section{Introduction}
In this paper, we give bijective proofs of some of  Amini's conjectures  \cite{amini} on equidistributions of Mahonian statistics over sets of pattern avoiding 
 permutations.

Let $S_n$ denote the set of all the permutations of $[n]=\{1,2,\cdots,n\}$.
Given a permutation $\pi=\pi_1 \pi_2 \cdots \pi_n \in S_n$,
let $\pi^c$ be the complement of $\pi$, which is given by
$\pi^c_i=n+1-\pi_i$. The reverse of $\pi$ is denoted $\pi^r$, where $\pi^r_i=\pi_{n+1-i}$.
A descent of $\pi$ is a position $i \in [n-1]$ such that $\pi_i > \pi_{i+1}$.
The descent set  of $\pi$
is given by
\[
D(\pi)=\{i\colon \pi_i>\pi_{i+1}\},\]
Define the set of the inversions of $\pi$ by
$$
Inv(\pi)=\{(i,j)\colon 1\leq i< j \leq n, \pi_i > \pi_j\}.
$$
Let $des(\pi)$and $inv(\pi)$
be the descent number and the inversion number of $\pi$, which are defined by
$des(\pi)=|D(\pi)|$ and $inv(\pi)=|Inv(\pi)|$, respectively.
The major index of  $\pi$, denoted $maj(\pi)$, is given by
\[maj(\pi)=\sum_{i \in  D(\pi)}  i.\]
A run of a permutation $\pi \in S_n$  is a maximal factor
(subsequence of consecutive elements) which is increasing.
For instance,  the permutation
$7356412$  has four runs, namely, $7$, $356$, $4$ and $12$.

Suppose that $st_1$ is a statistic on the object $Obj_1$ and $st_2$ is a statistic on the object $Obj_2$. If
\[
\sum_{\sigma \in Obj_1} q^{st_1(\sigma)}=\sum_{\sigma \in Obj_2} q^{st_2(\sigma)},
\]
we say that the statistic $st_1$ over $Obj_1$ is equidistributed with the statistic $st_2$ over $Obj_2$.

A statistic on $S_n$ is said to be Mahonian if
it is equidistributed with the statistic  $inv$ on $S_n$. It is well-known that
\[\sum_{\pi \in S_n} q^{inv(\pi)}=\sum_{\pi \in S_n} q^{maj(\pi)}=[n]_q !,\]
where $[n]_q=1+q+\cdots +q^{n-1}$ and $[n]_q!=[n]_q [n-1]_q \cdots [1]_q$. Thus, the major index $maj$ is a Mahonian statistic.

Given a permutation $\pi=\pi_1\pi_2\cdots \pi_n$ in $S_n$ and  a permutation $\sigma=\sigma_1 \sigma_2 \cdots \sigma_k \in S_k$, where $k\leq n$, we say
 that $\pi$ contains a pattern $\sigma$  if there exists a subsequence
 $\pi_{i_1}\pi_{i_2} \cdots \pi_{i_k}$ $(1 \leq i_1 < i_2 < \cdots <i_k \leq n)$ of $\pi$ that is order isomorphic to $\sigma$,
 in other words,  for  $1 \leq  l<m \leq k $,
 $\pi_{i_l}<\pi_{i_m}$ if and only if $\sigma_l <\sigma_m$. Otherwise, we say that $\pi$  avoids a pattern $\sigma$, or
 $\pi$ is $\sigma$-avoiding. For example, $74538126$ is 1234-avoiding, while it
  contains a pattern $3142$ corresponding to the subsequence $7486$.

In \cite{Babson}, Babson and Steingr\'{\i}msson introduced generalized permutation patterns, where
they write two adjacent letters
may or may not be separated by a dash.
Two adjacent letters without a dash in a pattern indicates that the corresponding letters in the permutation must be
adjacent. 
 For example, an occurrence of the generalized pattern $2$-$31$ in a permutation
$\pi=\pi_1\pi_2 \cdots \pi_n$ is a subsequence $\pi_i\pi_j\pi_{j+1}$ such that
$i<j$ and $\pi_{j+1}<\pi_i<\pi_j$.
Let $(\tau)\pi$ denote
 the number of occurrences of $\tau$ in $\pi$. Here, we see $(\tau)$ as a
  generalized pattern function.
For $\pi=4753162$, we have $(2$-$31)\pi=4$.

 Further, Babson and Steingr\'{\i}msson \cite{Babson} 
investigated the linear
combinations of generalized patterns. Some of them are listed below.
\begin{align*}
  maj= & \, (1-32)+(2-31)+(3-21)+(21), \\[3pt]
  stat= &\,(13-2)+(21-3)+(32-1)+(21),\\[3pt]
  mak= &\, (1-32)+(31-2)+(32-1)+(21).
\end{align*}
They
showed that almost all of
the Mahonian permutation statistics in the literature can be written as linear
combinations of generalized patterns. 
They conjectured that
other combinations of generalized patterns are still Mahonian. Some of the conjectures were proved and reproved
by Foata and Randrianarivony \cite{randrian}, Vajnovszki \cite{vaj}, etc. 

Particularly, Babson and Steingr\'{\i}msson  conjectured that the statistic $(des,stat)$ is Euler-Mahonian.
In 2001, D. Foata and  D. Zeilberger \cite{Foata} gave a proof of this conjecture  using
q-enumeration and generating functions and an almost completely automated proof via Maple packages ROTA and PERCY.

Given a permutation $\pi=\pi_1 \pi_2 \cdots \pi_n$, let $F(\pi)=\pi_1$ be the
first letter of $\pi$ and
\[adj(\pi)=|\{i \colon 1 \leq i \leq n \text{~and~} \pi_i- \pi_{i+1}=1\}|,\]
where $\pi_{n+1}=0$.
 Burstein \cite{Burstein} provided a
bijective proof of the following  refinement of Babson and Steingr\'{\i}msson's Conjecture  as follows.
\begin{thm}
Statistics $(adj, des, F, maj, stat)$ and $(adj, des, F, stat, maj)$
are equidistributed over $S_n$ for all $n$.
\end{thm}

In 2017, the author and Li \cite{chen} gave a new bijective proof of Babson and Steingr\'{\i}msson's Conjecture by analogy with Carlitz's insertion method \cite{Carlitz}.
While this bijection does not  preserve the statistic $adj$.

Recently, Amini \cite{amini} studied the distributions of  these Mahonian statistics defined by combinations of  generalized pattern functions at most three over sets of pattern avoiding permutations. 
He proved a host of equidistributions by adopting 
block decomposition, Dyck path and generating functions. 
Meanwhile, he gave some conjectures on it. We list some of them below. Note that Amini used the notation $bast$
instead of $stat$. 

\begin{conj}\label{mainconj}
The statistic $maj$ (resp. $stat$) over $S_n(213)$ is equally distributed
with the statistic $stat$ (resp. $maj$) over $S_n(231)$, namely,
\begin{align*}
  \sum_{\pi \in S_n(213)}q^{maj(\pi)}= & \sum_{\pi \in S_n(231)}q^{stat(\pi)} \\[5pt]
  \sum_{\pi \in S_n(231)}q^{maj(\pi)}= & \sum_{\pi \in S_n(213)}q^{stat(\pi)}.
\end{align*}
\end{conj}

\begin{conj}\label{conseconj1}
The statistic $maj$  over $S_n(132)$ (resp. $S_n(312)$) is equally distributed
with the statistic $stat$  over $S_n(213)$ (resp. $S_n(231)$). 
\end{conj}

\begin{conj}\label{conseconj2}
The statistic $mak$  over $S_n(132)$ (resp. $S_n(213)$, $S_n(231)$, $S_n(312)$ ) is equidistributed 
with the statistic $stat$  over $S_n(213)$ (resp. $S_n(231)$, $S_n(231)$, $S_n(213)$ ).
\end{conj}

In this paper, we will construct a bijection from $S_n(213)$
to $S_n(231)$, which settles Conjecture \ref{mainconj}. 
Then, Conjectures \ref{conseconj1} and Conjectures \ref{conseconj2} follows directly from Conjecture \ref{mainconj} and some known equidistributions 
given by Amini \cite{amini}.

\section{A bijective proof of Conjecture \ref{mainconj}}

In this section, we will give a bijective proof of the following theorem, which settles Conjecture \ref{mainconj}.
\begin{thm}\label{main}
Statistics $(F,maj,stat,des,ides)$ over $S_n(213)$ are
equidistributed with statistics $(F,stat,maj,des,ides)$
over $S_n(231)$ for all $n$.
\end{thm}


To construct the bijection which proves the above theorem,
we first define two maps $\phi$ and $\varphi$.

Let $S=\{s_1, s_2 ,\cdots, s_n\}$ with $s_1 < s_2 < \cdots <s_n$ and $T=\{t_1, t_2 , \ldots, t_n\}$ with $t_1 < t_2 <\cdots <t_n$. For a
permutation $\pi=\pi_1 \pi_2 \cdots \pi_n$ of the set $S$, assume that $\pi_i=s_1$ and $\pi=\pi' \pi_i \pi''$, where $\pi'$ and $\pi''$ may be empty. Let $T_{1}=\{t_1,\ldots, t_{n-i}\}$ and $T_{2}=\{t_{n-i+1},\ldots, t_{n-1}\}$.
When $i=1$\,($i=n$), we treat $T_2$\,($T_1$) as an empty set.
Define a map $\phi$ recursively by \[\phi(\pi,T)=\phi(\pi'',T_1) \,t_n\, \phi(\pi',T_2).\]
By the construction of $\phi$, it is easily checked that $\phi(\pi,T)$ is $231$-avoiding.
As an example, let $S=\{2,3,4,5,6\}$ and $T=\{1,3,5,7,9\}$.
It can be computed that
\begin{align*}
  \phi(63542,T) =& 9 \phi(6354,\{1,3,5,7\})\\
  = & 9 \phi(54,\{1,3\})7\phi(6,\{5\})\\
  =&93175.
\end{align*}
Given a permutation $\sigma=\sigma_1 \sigma_2 \cdots \sigma_n$ of the set $T$, assume that $\sigma_j=t_n$ and $\sigma=\sigma' \sigma_j \sigma''$, where $\sigma'$ and $\sigma''$ may be empty. Let $S_1=\{s_{j+1},\ldots, s_{n}\}$ and $S_2=\{s_{2},\ldots, s_{j}\}$.
When $j=1$\,($j=n$), we view $S_2$\,($S_1$) as an empty set. Define a map $\varphi$ recursively by \[\varphi(\sigma,S)=\phi(\sigma'',S_1) \,s_1\, \phi(\sigma',S_2).\]
Clearly, $\varphi(\sigma,S)$ is $213$-avoiding.
For example, it can be deduced that
\begin{align*}
  \varphi(93175,S) =& \varphi(3175,\{3,4,5,6\})\,2\\
  = & \varphi(5,\{6\})\, 3\,\varphi(31,\{4,5\})\,2\\
  =&63542.
\end{align*}
We notice that $\varphi(\phi(63542,T),S)=63542$,
while $63542$ is $213$-avoiding. In fact,
we have the following properties of maps $\phi$ and $\varphi$.
\begin{lem}\label{phivarphi}
Suppose $S$ and $T$ are both $n$-element sets. Given a
$213$-avoiding permutation $\pi$ of $S$, we have
$\varphi(\phi(\pi,T),S)=\pi$. Moreover, if $D(\pi)=\{j_1, j_2 ,\ldots, j_r\}$,  we have
\[D(\phi(\pi,T))
=\{n-j_r, \ldots,n-j_2 , n-j_1\}.\]
\end{lem}
\pf First, we proceed to prove $\varphi(\phi(\pi,T),S)=\pi$ by induction on $n$.
When $n=1$, it is trivial. Suppose it holds for
all $m$, where $m < n$, we wish to show that it holds for $n$.

Write $S=\{s_1, s_2 ,\cdots, s_n\}$ with $s_1 < s_2 < \cdots <s_n$ and $T=\{t_1, t_2 , \ldots, t_n\}$ with $t_1 < t_2 <\cdots <t_n$. Assume that $\pi=\pi' \pi_i \pi''$, where
$\pi_i=s_1$ and $\pi'$, $\pi''$ may be empty. Let $T_{1}=\{t_1,\ldots, t_{n-i}\}$ and $T_{2}=\{t_{n-i+1},\ldots, t_{n-1}\}$.
We have $\phi(\pi,T)=\phi(\pi'',T_1) \,t_n\, \phi(\pi',T_2)$.
It follows that
\[
\varphi(\phi(\pi,T),S) = \varphi(\phi(\pi',T_2),S') \,s_1 \,\varphi(\phi(\pi'',T_1),S''),
\]
where $S'=\{s_{n-i+2}, \ldots, s_n\}$ and $S''=\{s_2, \ldots, s_{n-i+1}\}$. Since $\pi$ is $213$-avoiding and $\pi_i=s_1$, we see that $\pi'$ is a $213$-avoiding permutation of $S'$ and $\pi''$ is a $213$-avoiding permutation of $S''$.
By the hypothesis of the induction, we deduce  that
$\varphi(\phi(\pi',T_2),S')=\pi'$ and $\varphi(\phi(\pi'',T_1),S'')=\pi''$. Hence,
$\varphi(\phi(\pi,T),S)=\pi' s_1 \pi''=\pi$.

In the following, we aim to prove that for $213$-avoiding permutation $\pi$ of $S$ with descent set $\{j_1, j_2 ,\ldots, j_r\}$,  we have
 \[D(\phi(\pi,T))=\{n-j_r, \ldots,n-j_2 , n-j_1\}.\]
 We proceed  by induction on $n$.

When $n=1$, we see that $S=\{s_1\}$ and $T=\{t_1\}$.
Notice that $\tau(s_1)=t_1$, $D(s_1)=\varnothing$ and
$D(t_1)=\varnothing$. Hence, this property holds for $n=1$.
Assume that it holds for
$m < n$. We wish to show that it holds for $n$.

 Recall that $\pi=\pi' \pi_i \pi''$, where $\pi_i=s_1$. If $i=1$, then $\pi'$ is empty and $D(\pi'')=\{j_1-1, \ldots, j_r-1\}$. By the assumption of this induction, we
  deduce that
  \begin{align*}
    D(\phi(\pi'',\{t_1,\ldots, t_{n-1}\}))=& \{n-1-{j_r-1}, \ldots, n-1-(j_1-1)\} \\
    =&\{n-j_r, \ldots, n-j_1 \}.
  \end{align*}
  Since $\phi(\pi,T)=\phi(\pi'',\{t_1,\ldots, t_{n-1}\})\,t_n$,
 it follows that \[D(\phi(\pi,T))=\{n-j_r, \ldots, n-j_1 \}.\]
  This completes the proof of this case.

   If $i>1$, we deduce that  $i-1 \in D(\pi)$. Let $j_d=i-1$. Then, we have $D(\pi')=\{j_1, \ldots,j_{d-1}\}$ and $D(\pi'')=\{j_{d+1}-i, \ldots,j_{r}-i\}$. It should be mentioned that if $d-1<1$\,($d+1>r$), we treat $D(\pi')$\,($D(\pi'')$) as an empty set.
   Recall that  \[\phi(\pi,T)=\phi(\pi'',T_1)\, n \, \phi(\pi',T_2).\]
   We deduce that $n-i+1 \in D(\phi(\pi,T))$, namely,
   $n-j_d \in D(\phi(\pi,T))$. Moreover, bazy the assumption of this
   induction, we  see that
   \begin{align*}
     D(\phi(\pi'',T_1))= & \{n-i-(j_r-i),\ldots,n-i-(j_{d+1}-i)\} \\
     = & \{n-j_r, \ldots, n-j_{d+1}\}
   \end{align*}
 and
 \begin{align*}
     D(\phi(\pi',T_2))= & \{i-1-j_{d-1},\ldots,i-1-j_{1}\}.
   \end{align*}
  It follows that
\begin{align*}
  D(\phi(\pi,T))= & \{n-j_r, \ldots, n-j_{d+1},
  n-j_d, n-i+1+(i-1-j_{d-1}), \ldots, n-i+1+i-1-j_{1}\} \\
  = & \{n-j_r, \ldots,n-j_2 , n-j_1\}.
\end{align*}
This completes the proof.\qed

Clearly, the map $\phi$ keeps the statistic $des$ of the permutations in its preimage and image. Moreover, $\phi$
also keeps $ides$. To show this, we need the following lemmas.

\begin{lem}\label{213des}
For $\pi=\pi_1 \pi_2 \cdots \pi_n \in S_n(213)$, we have $des(\pi)=ides(\pi)$.
\end{lem}
\pf Assume that $D(\pi)=\{j_1,j_2,\ldots,j_r\}$. Then, the runs of $\pi$ are $\pi_1 \cdots \pi_{j_1}$,\, $\pi_{j_1+1} \cdots \pi_{j_2}$, \ldots, $\pi_{j_{r-1}+1} \cdots \pi_{j_r}$,\, $\pi_{j_{r}+1} \cdots \pi_n$.
In the following, we proceed to show that  $D(\pi^{-1})=\{\pi_n,\pi_{j_r}, \ldots, \pi_{j_2}\}$.

Notice that $i \in D(\pi^{-1})$ if and only if
$i$ is located on the right of $i+1$ in $\pi$.
 For $2 \leq k \leq r$,
we claim that $\pi_{j_k}$ is on the right of $\pi_{j_k}+1$.
If not, $\pi_{j_k}\pi_{j_k+1}\,(\pi_{j_k}+1)$ will form
a $213$-pattern, which contradicts the fact that $\pi$ is $213$-avoiding. The claim is verified. It follows that $\pi_{j_k} \in D(\pi^{-1})$. Notice that  $\pi_n \in D(\pi^{-1})$. Thus, we deduce that $\{\pi_n,\pi_{j_r}, \ldots, \pi_{j_2}\} \subseteq D(\pi^{-1})$.

For ${j_k} < l < j_{k+1}$ ($1 \leq k \leq r-1$) or
$ j_r < l<n$, we claim that $\pi_l \notin D(\pi^{-1})$.
Namely, $\pi_l$ is on the left of
$\pi_l+1$. Otherwise, $(\pi_l+1) \pi_l \pi_{l+1}$ will
form  a $213$-pattern.  Combining with the fact that $\pi_l$ ($1 \leq l \leq j_1$) is on the left of $\pi_l+1$,
we deduce that $D(\pi^{-1})=\{\pi_n,\pi_{j_r}, \ldots, \pi_{j_2}\}$. This completes the proof. \qed

By the fact that $(\pi^c)^{-1}=(\pi^{-1})^r$,
the   property  below follows directly from Lemma \ref{213des}.
\begin{cor}\label{231des}
For $\pi \in S_n(231)$, we have $des(\pi)=ides(\pi)$.
\end{cor}

As a direct deduction of Lemma \ref{phivarphi}, \ref{213des} and Corollary \ref{231des}, we have the following properties.
\begin{cor}\label{desides}
Let  $S$ and $T$ be $n$-element sets. For a $213$-avoiding
permutation $\pi$ of $S$, we have
\begin{align*}
  des(\pi) = & \,des(\phi(\pi,T)), \\[3pt]
  ides(\pi)= & \,ides(\phi(\pi,T)).
\end{align*}
\end{cor}

Now, we  are ready to give a description of the desired map $\alpha$ from $S_n(213)$ to $S_n(231)$, which turns out to give a bijective proof of Theorem \ref{main}. For a permutation $\pi \in S_n(213)$,
suppose $F(\pi)=k$.
We write $\pi=k\pi' \pi''$, where $\pi'$ is a permutation
of $\{k+1,\ldots, n\}$ and $\pi''$ is a permutation
of $\{1,\ldots,k-1\}$.
Both $\pi'$ and $\pi''$ are allowed to be empty. Let \[\alpha(\pi)=k\,\phi(\pi'',[k-1])\,\phi(\pi',\{k+1,\ldots,n\}).\]
For instance, given a permutation $25678341 \in S_8(213)$, we deduce that
\begin{align*}
  \alpha(25678341)= & 2\,\phi(1,\{1\})\,\phi(567834,\{3,4,5,6,7,8\}) \\
  =& 21384567.
\end{align*}

\begin{thm}\label{bijection}
$\alpha$ is a bijection from $S_n(213)$ to $S_n(231)$.
\end{thm}
\pf Given $\pi \in S_n(213)$, write $\pi=k \pi' \pi''$ as above.
By the construction of $\phi$, we see that
$\phi(\pi'',[k-1])$ is a $231$-avoiding permutation of $[k-1]$, while $\phi(\pi',\{k+1,\ldots,n\})$ is a $231$-avoiding permutation of $\{k+1,\ldots,n\}$.
Following from the fact that $\alpha(\pi)=k\,\phi(\pi'',[k-1])\,\phi(\pi',\{k+1,\ldots,n\})$,
we deduce that $\alpha(\pi)$ is $231$-avoiding.
Hence, $\alpha$ is a map from $S_n(213)$ to $S_n(231)$.

To prove that $\alpha$ is a bijection, we construct its inverse. For a permutation $\sigma$ in $S_n(231)$,
assume that $F(\sigma)=f$ and $\sigma=f \sigma' \sigma''$,
where $\sigma'$ is a $231$-avoiding permutation of
$\{1,\ldots,f-1\}$ and
$\sigma''$ is a $231$-avoiding permutation of $\{f+1, \ldots, n\}$. Both $\sigma'$ and $\sigma''$ are allowed to be empty.
Define $\beta(\sigma)$ by
\[\beta(\sigma)=f \, \varphi(\sigma'',\{f+1, \ldots, n\})\, \varphi(\sigma',[f-1]).\]
We claim that $\beta$ is the inverse map of $\alpha$.
It suffices to prove that $\beta(\alpha(\pi))=\pi$.
Notice that
\begin{align*}
  \beta(\alpha(\pi))= & \beta(k\,\phi(\pi'',[k-1])\,\phi(\pi',\{k+1,\ldots,n\})) \\
  = & k \, \varphi(\phi(\pi',\{k+1,\ldots,n\}),\{k+1,\ldots,n\})\, \varphi(\phi(\pi'',[k-1]),[k-1])\\
  =& k \pi' \pi''\\
  =& \pi.
\end{align*}
Hence, the claim is verified and we compete the proof.\qed

As an immediate consequence, we see that the statistic $F$ is equidistributed over $S_n(213)$ and $S_n(231)$. To give a
proof of Theorem \ref{main}, we need the following lemmas,
the first of which was given by Burstein \cite{Burstein}.

\begin{lem}\label{Bur}
For $\pi \in S_n$ , we have \[maj(\pi) + stat(\pi) = (n + 1) des (\pi) - (F(\pi) - 1).\]
\end{lem}
\begin{lem}\label{equi}
Given $\pi \in S_n(213)$, we have
\begin{align*}
  des(\pi)=& \, des(\alpha(\pi)),\\[2pt]
  ides(\pi)=& \, ides(\alpha(\pi)),\\[2pt]
  maj(\pi) =& \, stat(\alpha(\pi)),\\[2pt]
   stat(\pi) =& \, maj(\alpha(\pi)).
\end{align*}
\end{lem}

\pf Assume that $\pi=k \pi'\pi''$, where $\pi'$ is a permutation of $\{k+1, \ldots, n\}$ and $\pi''$ is a permutation of $\{1, \ldots, k-1\}$.
By definition, we have
\[\alpha(\pi)=k\,\phi(\pi'',[k-1])\,\phi(\pi',\{k+1,\ldots,n\}).\]
It follows that
\[des(\alpha(\pi)) = 1+ des(\phi(\pi'',[k-1]))+ des(\phi(\pi',\{k+1,\ldots,n\})).\]
By Corollary \ref{desides}, we see that
$des(\phi(\pi'',[k-1]))=des(\pi'')$ and $des(\phi(\pi',\{k+1,\ldots,n\}))=des(\pi')$.
Hence, we have  $des(\alpha(\pi))=1+des(\pi'')+des(\pi')=des(\pi)$.

By Lemma   \ref{213des}, we see that
$des(\pi)=ides(\pi)$. By Corollary \ref{231des}, we
have $des(\alpha(\pi))=ides(\alpha(\pi))$.
Thus, \,$ides(\pi)= ides(\alpha(\pi))$ follows.

In the following, we proceed to prove that $maj(\pi) = stat(\alpha(\pi))$. We consider two cases.
When $n=1$, $\pi=1 \pi'$. Assume that
$D(\pi)=\{j_1, j_2, \cdots, j_r\}$. Then, we have
$D(\pi')=\{j_1-1, j_2-1, \cdots, j_r-1\}.$
By Lemma \ref{phivarphi}, we see that
$D(\phi(\pi', \{2, \ldots, n\}))=\{n-j_r, \ldots, n-j_1\}$.
Hence, the descent set of $\alpha(\pi)=1 \, \phi(\pi',\{2, \ldots, n\})$ is $\{n+1-j_r, \ldots, n+1-j_1\}$. It follows that
\begin{align*}
 maj(\alpha(\pi))= & \, (n+1-j_r)+ \cdots+(n+1-j_1) \\[3pt]
  = & \, (n+1)des(\alpha(\pi))-maj(\pi)
\end{align*}
By Lemma \ref{Bur} and the fact $F(\pi)=1$, we deduce that
\begin{align*}
  stat(\alpha(\pi))= &\, (n+1)des(\alpha(\pi))- maj(\alpha(\pi))\\
  =& \, maj(\pi).
\end{align*}
This completes the proof of the case $n=1$.

When $n \neq 1$, suppose that the descent set of
$\pi$ is $\{j_1, \ldots, j_{d-1},\, j_d,\, j_{d+1}, \ldots, j_r\}$, where $j_d=n-k+1$. Thus, we deduce that
$D(\pi')=\{j_1-1, j_2-1, \ldots, j_{d-1}-1\}$
and $D(\pi'')=\{j_{d+1}-(n-k+1), \ldots, j_{r}-(n-k+1)\}$.
By Lemma \ref{phivarphi}, we have
\[D(\phi(\pi',\{k+1, \ldots, n\}))=\{n-k-j_{d-1}+1, \ldots, n-k-j_1+1\},\]
and
\begin{align*}
  D(\phi(\pi'',[k-1]))= & \,\{k-1-(j_{r}-(n-k+1)),\ldots, k-1-(j_{d+1}-(n-k+1))\}     \\
  =& \,\{n-j_r, \ldots, n-j_{d+1}\}.
\end{align*}
By the definition of the map $\alpha$, it follows that
\begin{align*}
  D(\alpha(\pi)) =& \, \{1,\, n+1-j_r, \ldots, n+1-j_{d+1}, \,k+(n-k-j_{d-1}+1), \ldots, k+(n-k-j_{1}+1)\} \\[3pt]
 = & \, \{1, \,n+1-j_r, \ldots, n+1-j_{d+1}, \,n+1-j_{d-1}, \ldots, n+1-j_1\}.
\end{align*}
Hence, we  deduce that
 \begin{align*}
         maj(\alpha(\pi))= & 1+(n+1-j_r)+\cdots+ (n+1-j_{d+1})+ \cdots+ (n+1-j_{d-1})+ \cdots+(n+1-j_1)  \\[3pt]
         = & (n+1) (des(\alpha(\pi))-1) +1- j_1 - \cdots -j_{d-1}-j_{d+1}-\cdots-j_r
       \end{align*}
Since $maj(\pi)=j_1 + \cdots +j_{d-1}+j_{d+1}+\cdots+j_r +(n-k+1)$ and $des(\alpha(\pi))=des(\pi)$, we have
\begin{equation}\label{c1}
  maj(\alpha(\pi))+maj(\pi)=(n+1)des(\pi)-(k-1).
\end{equation}
Moreover, by Lemma \ref{Bur} we have
\begin{equation}\label{c2}
  maj(\alpha(\pi))+stat(\alpha(\pi))=(n+1)des(\pi)-(k-1).
\end{equation}
Combining (\ref{c1}) and (\ref{c2}), we deduce that  $maj(\pi)=stat(\alpha(\pi))$.
This completes the proof of the case $n \neq 1$.

Notice that by Lemma \ref{Bur}, we also have
\begin{equation}\label{c3}
  maj(\pi)+stat(\pi)=(n+1)des(\pi)-(k-1).
\end{equation}
Combining relations (\ref{c1}) and  (\ref{c3}), we deduce that
\[stat(\pi)=maj(\alpha(\pi)).\]
We compete the proof. \qed

By Theorem \ref{bijection} and Lemma \ref{equi}, we give a bijective proof of Theorem \ref{main}.
As mentioned in the introduction, there are
two involutions on $S_n$ which map the
statistic $(maj,stat)$ to $(stat, maj)$ given
by Burstein \cite{Burstein} and Chen and Li \cite{chen}, respectively. We note that
$\alpha$ is not the restriction of either of them.

At the end of this paper, we give a consequence of a fact in the proof of Lemma \ref{equi}.
\begin{cor} For $n \geq 1$, we have
\begin{equation}\label{Mn}
  M_n(q,t)=M_{n-1}(q,t)+\sum_{k=1}^{n-1}q^{k}tM_{k-1}(q,qt)M_{n-k}(q,q^kt),
\end{equation}
where \[M_n(q,t)=\sum_{\pi \in S_n(213)}q^{maj(\pi)}t^{des(\pi)}.\]
\end{cor}
\pf Given $\pi \in S_n(213)$, use the notations $\pi=k\pi'\pi''$ of the previous proof. We have
\begin{equation}\label{recur1}
des(\pi)= \left\{
  \begin{array}{ll}
des(\pi'), & \mbox{$k=1$}, \\[3pt]
1+ des(\pi')+des(\pi''), & \mbox{$k \neq 1$.}
 \end{array} \right.
\end{equation}
Moreover, we have 
\begin{equation}\label{recur2}
  maj(\pi)=maj(\pi')+maj(\pi'')+des(\pi')+des(\pi'')(n-k+1)+n-k+1.
\end{equation}
Then, the recurrence relation (\ref{Mn}) comes directly from relations
(\ref{recur1}) and (\ref{recur2}). \qed

It should be mentioned that similarly as given in \cite{sagan},
 by another block decomposition of $\pi$
 as $\pi=\pi' 1 \pi''$, we may obtain the following 
 recurrence relations of $M_n(q,t)$, which is different
 from (\ref{Mn}),
 \[M_n(q,t)=M_{n-1}(q,t)+\sum_{k=1}^{n-1}q^{k}tM_{k}(q,t)M_{n-k-1}(q,q^{k+1}t).\]


\end{document}